%marko.tex:  Counting Clean Words According to the Number of Their Clean Neighbors
%%a Plain TeX file  Doron Zeilberger (6 pages)

%begin macros

\baselineskip=14pt
\parskip=10pt

\font\eightrm=cmr8 

\magnification=\magstephalf

\def\C{{\cal C}}

\def\1{{\overline{1}}}
\def\2{{\overline{2}}}
\parindent=0pt
\overfullrule=0in

\def\frac#1#2{{#1 \over #2}}
%\headline={\rm  \ifodd\pageno  \RightHead  \else  \LeftHead  \fi}
%\def\RightHead{\centerline{
%Title
%}}
%\def\LeftHead{ \centerline{Doron Zeilberger}}
%end macros
\centerline
{\bf 
Counting Clean Words According to the Number of Their Clean Neighbors
}
\bigskip
\centerline
{\it Shalosh B. EKHAD and Doron ZEILBERGER}
\bigskip
\qquad {\it In fond memory of Marko Petkov${\breve s}$ek (1955-2023), a great summer and enumerator}

{\bf Preface}

Our good friend and collaborate, Marko Petkov${\breve s}$ek ([PWZ]), passed away on March 23, 2023, 
and we already wrote a eulogy [Z], and donated to the Online Encyclopedia of Integer Sequences in his memory (See {\tt https://oeisf.org/donate} and search for Petkovsek).
However we believe that we can do more than that to commemorate Marko. We looked through his list of publications, and found the delightful
article [KMP] by Marko, joint with Sandi Klav${\breve z}$ar and Michel Mollard, and realized that the beautiful methodology that they used to solve
one very {\it specific} enumeration problem is applicable to a wide class of enumeration problems of the same flavor. More important, since
Marko was such an authority in {\it symbolic} computation, we decided to {\bf implement} the method, and wrote a Maple package

{\tt https://sites.math.rutgers.edu/\~{}zeilberg/tokhniot/Marko.txt} \quad,

that can very fast answer these kind of questions. In particular as we will soon see, Theorem 1.1 of [KMP] can be gotten (in its equivalent
form in terms of generating functions stated as $f(x,y)$ on top of p. 1321) by typing

{\tt WtEs( $\{0$,1 $\}$,$\{$[1,1]$\}$,y,x,3);} \quad .

Our Maple package, {\tt Marko.txt}, gives, in $0.057$ seconds, the answer
$$
-\frac{x^{2} y^{2}-x^{2} y -x y -1}{x^{3} y^{2}-x^{3} y -x^{2} y -x y +1} \quad.
$$

{\bf The Problem Treated so Nicely by Klav$ {\bf {\breve z}}$ar, Mollard, and Petkov${\breve s}$ek}

There are $2^n$ vertices in the $n$-dimensional unit cube $\{0,1\}^n$ and every such vertex has {\bf exactly} $n$ neighbors
(i.e. vertices with Hamming distance $1$ from it). The {\bf Fibonacci lattice} consists of those vertices whose
$01$ vector {\bf avoids} two consecutive $1$s, in other words of words in the alphabet $\{0,1\}$ avoiding
as a {\bf consecutive subword} the two-letter word $11$. Such words are called {\bf Fibonacci words}, and
there are, not surprisingly, $F_{n+2}$ of them (why?).

Each such word has $n$ neighbors, but some of them are not Fibonacci words. The question answered so elegantly in [KMP] was:

For any given $n$ and $k$, How many Fibonacci words of length $n$ are there that have exactly $k$ Fibonacci neighbors?
Calling this number $f_{n,k}$, [KMP] derived an explicit expression for it, that is equivalent to the
generating function (that they also derived)
$$
f(x,y) \, = \, \sum_{n,k \geq 0} \, f_{n,k} \, x^n\, y^k \, = \,
-\frac{x^{2} y^{2}-x^{2} y -x y -1}{x^{3} y^{2}-x^{3} y -x^{2} y -x y +1} \quad.
$$

They also considered the analogous problem for {\bf Lucas words} that consists of Fibonacci words where the first and last letter can't both be $1$.
This problem is also amenable to far-reaching generalization, but will not be handled here.

{\bf The general Problem}

{\bf Input}: 

$\bullet$ A finite alphabet $A$ (In the [KMP] case $A=\{0,1\}$).

$\bullet$ A finite set of words $M$, (of the same length) in the alphabet $A$. (In the [KMP] case $M$ is the singleton set $\{11\}$).

{\bf Definition}: A word in the alphabet $A$ is called {\bf clean} if it does not have, as {\it consecutive substring}, any of the members of $M$.

In other words writing $w=w_1 \dots w_n$, a word is {\bf dirty} if there exists an $i$ such that $w_i w_{i+1} \dots w_{i+k-1} \in M$.
For example if $A=\{1,2,3\}$ and $M=\{123,213\}$, then $12212312$ is dirty while $111222333$ is clean.

{\eightrm To get the set of clean words of length {\tt n} in the alphabet {\tt A} and set of `mistakes' {\tt M}, type, in {\tt Marko.txt},

{\tt CleanWords(A,M,n);} \quad .

For example, to get the Fibonacci words of length $3$ type:

{\tt CleanWords($\{$ 0,1 $\}$,$\{$ [1,1] $\}$ , 3);}, getting:

$\{$[0, 0, 0], [0, 0, 1], [0, 1, 0], [1, 0, 0], [1, 0, 1]$\}$ \quad .

}

The problem of the straight enumeration of clean words is handled very efficiently via the Goulden-Jackson cluster algorithm [NZ], but it is not suitable for the present problem of
{\it weighted} enumeration.

{\bf Definition}: Two words of the same length in the alphabet $A$ are neighbors if their Hamming distance is $1$, in other words,
$u=u_1 \dots u_n$ and $v=v_1 \dots v_n$ are neighbors if there exists a location $r$ such $u_i=v_i$ if $i \neq r$ and $u_r \neq v_r$.

For example if $A=\{1,2,3\}$, the set of neighbors of $111$  is
$$
\{ 211,311,121,131,112,113   \} \quad .
$$
Obviously every word of length $n$ in the alphabet $A$ has $n\cdot (|A|-1)$ neighbors.

However, if $w$ is a clean word, some of its neighbors may be dirty, so if there is one {\it typo}, it can become dirty, and that would be embarrassing
(Oops, {\it embarrassing} is already dirty). While the word, {\it duckling} is clean, not all its
neighbors are clean.

{\eightrm

To see the number of clean neighbors of a word {\tt w} in the alphabet {\tt A} and set of mistakes {\tt M}, type

{\tt NCN(w,A,M); }
}

{\bf Output}: Having fixed the (finite) alphabet $A$, and the finite set of forbidden substrings $M$ (all of the same length),
let $f_{n,k}$  be the number of clean words in the alphabet $A$ of length $n$ having $k$ clean neighbors. Compute the
bi-variate generating function
$$
f(x,y) := \sum_{n,k \geq 0} f_{n,k} x^n\, y^k \quad .
$$
It would follow from the algorithm (inspired by the methodology of [KMP], but vastly generalized) that this is always a {\bf rational function} of $x$ and $y$.

{\eightrm This is implemented in procedure

{\tt WtEs(A,M,y,x,MaxK)},

where {\tt MaxK} is a `maximum complexity parameter'.
See the beginning of this article for the case treated in [KMP]. For a more complicated example, where a word is clean if it avoids the
substrings $000$ and $111$, type

{\tt WtEs( $\{$ 0,1 $\}$, $\{$ [1,1,1],[0,0,0] $\}$,y,x,5);}

getting, immediately:

$$
\frac{2 x^{5} y^{4}-4 x^{5} y^{3}+2 x^{5} y^{2}-2 x^{4} y^{3}+4 x^{4} y^{2}-2 x^{4} y -y^{2} x^{3}+2 x^{3} y -4 x^{2} y^{2}-x^{3}+2 x^{2} y +x^{2}-2 x y +x -1}{y^{2} x^{3}-x^{3}+x^{2}+x -1} \, .
$$

}

If you want to keep track of the individual letters, rather than just the length, use the more general procedure

{\tt WtEg(A,M,x,y,t,MaxK)}.

{\bf Reverse-Engineering the beautiful Klav${\bf {\breve z}}$ar-Mollard-Petkov${\bf {\breve s}}$ek Proof and Vastly Generalizing It}

In fact, the authors of [KMP] proved their results in two ways, and only the second way used {\it generatingfunctionology}. Even that part
argued directly in terms of the (double) sequence $f_{n,k}$ itself, and only at the {\bf end of the day}, took the
(bi-variate) generating function.

A more efficient, and {\it streamlined}, approach is to forgo the actual bi-sequence and operate {\bf directly} with weight-enumerators.
Let $\C(A,M)$ be the (`infinite') set of words in the alphabet $A$, avoiding, as consecutive substrings, the members of $M$, and for
each word $w$ in $\C(A,M)$, define the {\bf weight}, $Weight(w)$ by
$$
Weight(w)=x^{length(w)}\, y^{NCN(w)} \quad.
$$

For example, for the original case of $A=\{0,1\}$ and $M=\{11\}$,

$$
Weight(10101)\, = \, x^5\,y^3 \quad .
$$

We are interested in the {\bf weight-enumerator}
$$
f(x,y):= Weight(\C(A,M))= \sum_{w \in \C(A,M)} Weight(w) \quad.
$$
Once you have it, and you are interested in a specific $f_{n,k}$, all you need is to take a Taylor expansion about $(0,0)$ and extract the coefficient of $x^ny^k$.

Let $\C(A,M)^{(i)}$ be the subset of $\C(A,M)$ of words of length $i$, and pick a positive integer $k$. For any word $v \in \C(A,M)^{(k)}$, let
$\C_v(A,M)$ be the set of words in $\C(A,M)$ of length $\geq k$ that start with $v$. Obviously
$$
\C(A,M)= \bigcup_{i=0}^{k-1} \C(A,M)^{(i)} \, \,\cup \, \,\bigcup_{v \in \C(A,M)^{(k)}} \C(A,M)_v \quad .
$$

We can decompose $\C(A,M)_v$ as follows
$$
\C(A,M)_v \, = \, \bigcup_{a \in A} \, \C(A,M)_{va} \quad,
$$
where, of course $\C(A,M)_{va}$ is empty if appending the letter $a$  turns the clean $v$ into a dirty word.
Now, writing $v=v_1 \dots v_k$, and for $a \in A$  the computer verifies whether the difference
$$
NCN(v_1 \dots v_k a w)- NCN(v_2 \dots v_k a w) \quad
$$
is always the same, for any $v_1 \dots v_k a w \in \C(A_M)_{va}$. The way we implemented it is to test it for sufficiently long words, and then in {\it retrospect} have the
computer check it `logically', by looking the at the difference in the number of clean neighbors that happens by deleting the first letter $v_1$.
Let's call this constant quantity $\alpha(v,a)$.

It follows that we have a system of $|\C(A,M)^{(k)}|$  equations with $|\C(A,M)^{(k)}|$ unknowns.

$$
Weight(\C(A,M)_v) \, = \,
\sum_{{{a \in A}\atop{ va \in \C(A,M)}}}
x y^{\alpha(v,a)} Weight(\C(A,M)_{v_2 \dots v_{k-1} a}) \quad.
$$

After the computer algebra system (Maple in our case) automatically found all the $\alpha(v,a)$, and set up the system of equations, we kindly asked it to {\tt solve} it, getting
certain rational functions of $x$ and $y$. {\bf Finally}, our object of desire, $f(x,y)$, is given by
$$
Weight(\C(A,M)) \,= \,
\sum_{i=0}^{k-1} Weight(\C(A,M)^{(i)}) \, + \,
\sum_{v \in \C(A,M)^{(k)}} Weight(\C(A,M)_v) \quad.
$$

{\eightrm This is implemented in procedure {\tt WtEs(A,M,y,x,MaxK)}. }

If you also want to keep track of the individual letters, having the variable $t$ take care of the length, the equations are

$$
Weight(\C(A,M)_v) \, = \,
\sum_{{{a \in A}\atop{ va \in \C(A,M)}}}
x_{v_1} t y^{\alpha(v,a)} Weight(\C(A,M)_{v_2 \dots v_{k-1} a}) \quad.
$$

{\eightrm This is implemented in procedure {\tt WtEg(A,M,x,y,t,MaxK)}. }

{\bf Sample output}

$\bullet$ If you want to see the bi-variate generating functions for words in the alphabet $\{0,1\}$, avoiding $i$ consecutive occurrences of $1$, for $2 \leq i \leq 6$, see

{\tt https://sites.math.rutgers.edu/\~{}zeilberg/tokhniot/oMarko1.txt} \quad .

Note that the original case was $i=2$.

$\bullet$ If you want to see the bi-variate generating functions for words in the alphabet $\{0,1\}$, avoiding $i$ consecutive occurrences of $1$, 
{\bf and} $i$ consecutive occurrences of $0$,  for $3\leq i \leq 6$, see

{\tt https://sites.math.rutgers.edu/\~{}zeilberg/tokhniot/oMarko2.txt} \quad .

$\bullet$ If you want to see all such generating functions (still with BINARY words) for all possible SINGLE patterns of length 3,4,5 (up to symmetry), look at:

{\tt https://sites.math.rutgers.edu/\~{}zeilberg/tokhniot/oMarko3.txt} \quad .

The front of this article contains numerous other output files, but you dear reader, can generate much more!

{\bf Conclusion}

The value of the article [KMP], that inspired the present article, is not so much with the actual result, that in hindsight, thanks to our Maple package, is trivial, but in 
the human-generated ideas and {\bf methodology} that enabled one of us to generalize it to a much more general framework.

{\bf References}

[KMP] Sandi Klav${\breve z}$ar, Michel Mollard and  Marko Petkov${\breve s}$ek, {\it The degree sequence of Fibonacci and Lucas cubes}, Discrete Math. {\bf 311} (2011) 1310-1322. \hfill\break
{\tt https://www-fourier.ujf-grenoble.fr/\~{}mollard/soumis/DegSeqSubmit.pdf}

[NZ] John Noonan and Doron Zeilberger, {\it The Goulden-Jackson Cluster Method: Extensions, Applications, and Implementations},
J. Difference Eq. Appl. {\bf 5} (1999), 355-377.\hfill\break
{\tt https://sites.math.rutgers.edu/\~{}zeilberg/mamarim/mamarimhtml/gj.html} \quad .

[PWZ] Marko Petkov${\breve s}$ek, Herbert S. Wilf, and Doron Zeilberger, {\it ``A=B''}, A.K. Peters, 1996.\hfill\break
{\tt https://www2.math.upenn.edu/\~{}wilf/AeqB.pdf} \quad .

[Z] Doron Zeilberger, {\it  Marko Petko${\breve s}$v${\breve s}$ek (1955-2023), My A=B Mate }, Personal Journal of Shalosh B. Ekhad and Doron Zeilberger, April 14, 2023. \hfill\break
{\tt https://sites.math.rutgers.edu/\~{}zeilberg/mamarim/mamarimhtml/mpm.html} \quad .

\bigskip
\hrule
\bigskip

Shalosh B. Ekhad, c/o D. Zeilberger, Department of Mathematics, Rutgers University (New Brunswick), Hill Center-Busch Campus, 110 Frelinghuysen
Rd., Piscataway, NJ 08854-8019, USA. \hfill\break
Email: {\tt ShaloshBEkhad at gmail  dot com}   \quad .
\bigskip

Doron Zeilberger, Department of Mathematics, Rutgers University (New Brunswick), Hill Center-Busch Campus, 110 Frelinghuysen
Rd., Piscataway, NJ 08854-8019, USA. \hfill\break
Email: {\tt DoronZeil at gmail  dot com}   \quad .
\bigskip

April 21, 2023.

\end